\documentclass[10pt]{amsart}

\usepackage{amsmath,graphics}
\usepackage{amsfonts,amssymb}
\usepackage{amsmath,amscd}

\theoremstyle{plain}
\newtheorem*{theorem*}{Theorem}
\newtheorem*{lemma*} {Lemma}
\newtheorem*{corollary*} {Corollary}
\newtheorem*{proposition*} {Proposition}
\newtheorem{theorem}{Theorem}[section]
\newtheorem{lemma}[theorem]{Lemma}

\newtheorem{proposition}[theorem]{Proposition}

\theoremstyle{remark}

\newtheorem*{definition}{Definition}

\theoremstyle{definition}

\def\part{\partial}

\def\bp{\begin{pmatrix}}

\def\ep{\end{pmatrix}}
\def\bn{\begin{enumerate}}

\def\en{\end{enumerate}}
\def\ba{\begin{array}}
\def\ea{\end{array}}

\def\fr12{\frac{1}{2}}

\def\cmtbf#1{} \def\cmt#1{}

\begin{document}

\title{Local constancy of dimensions of hecke eigenspaces of automorphic forms}
\author{Aftab Pande}
\date{\today}
\begin{abstract}

We use a method of Buzzard to study $p$-adic families of
Hilbert modular forms and modular forms over imaginary
   quadratic fields. In the case of Hilbert modular forms, we get local constancy of dimensions of
   spaces of fixed slope and varying weight. For imaginary quadratic
   fields we obtain bounds independent of the weight on the dimensions
   of such spaces.

 \end{abstract}
\maketitle

\section{Introduction}

In this paper we explore $p$-adic variations of automorphic forms.
 Serre \cite{Se1} first presented the notion of a $p$-adic
analytic family of modular eigenforms using  $p$-adic
Eisenstein series, and this provided the first application of $p$-adic families.
 Hida then showed the first example of families of cuspidal eigenforms. His results
(\cite{Hi3}, \cite{Hi4}), were limited to the case of ordinary
modular forms, but proved instrumental in many number-theoretical
applications.

There was a wait of about a decade before non-ordinary families were
constructed. Using a rigid-analytic method of overconvergent modular
forms (based on earlier work of Katz \cite{Ka}) , Coleman proved the
existence of many families. He showed that almost
 every overconvergent
eigenform of finite slope lives in a $p$-adic family. The slope of
an eigenform is the $p$-adic valuation of its $U_p$-eigenvalue, and
having finite slope is a vast generalization of being ordinary,
i.e., to have a $U_p$-eigenvalue which is a $p$-adic unit. Coleman
also showed that overconvergent modular forms of small slope are
classical, which showed the existence of $p$-adic families of
classical modular forms. Coleman's work was motivated by, and
answered, a variety of questions and conjectures that Gouvea and
Mazur \cite{GM} had made based on ample numerical evidence. Coleman and Mazur \cite{CM} organized Coleman's results
(and more) in the form of a geometric object which was called the
eigencurve. It is a rigid-analytic curve whose points correspond to
normalized finite-slope $p$-adic overconvergent modular eigenforms
of a fixed tame level $N$.

One of the questions which remained was how big the radius $r$ of
the disc corresponding to a family could be. In \cite{GM}, Gouvea
and Mazur made some precise conjectures based on a lot of numerical
computations. The exact conjecture was disproved by Buzzard and
Calegari but Wan \cite{W} showed, using Coleman's theory of rigid-analytic
methods, that an eigenform $f$ of slope $s$ should live in a family
of eigenforms with radius $p^{-t}$, where $t = O(s^{2})$. Using
fairly elementary methods of group cohomology Buzzard \cite{B1}
found explicit bounds for the number of forms of slope $\alpha$,
weight $k$ and level $Np$, independent of the weight $k$. In his
unpublished paper \cite{B2}, he showed that forms have some kind of
$p$-adic continuity and gets results similar to Wan. We use these methods in the case of Hilbert modular forms and modular
forms over imaginary quadratic fields .

Let $F$ be a totally real field, $[F:\mathbb{Q}] = d$, where $d$ is
even. Let $\mathbf{k} = (k_1, k_2,....,k_d)$. Then there is a notion of
Hilbert cusp forms of weight $\mathbf{k} $ and level $\mathfrak{n}$ (an ideal of $O_F$).
By Jacquet-Langlands we get a relation between this space and
modular forms over a totally definite quaternion algebra $D$. Let us
call that space $\mathcal{S}_{ \mathbf{k}}^D(U) $, where $U$ is a compact open
subgroup of $D_f^{*}$, the adelisation of $D$.

Let $p$ be a fixed rational prime inert in $F$ and denote by
 $T_p$ the Hecke operator on the space of automorphic forms $\mathcal{S}_{\overrightarrow{k}}^D(U) $. We give
this space an integral structure for a ring $R$, where
$O_{F,p} \subseteq R$ and call it $\mathcal{S}_{ \mathbf{k}}^D(U,R) $.
There is a description of this space in terms of
$H^0(\Gamma^i(U),R)$, where the $\Gamma^i(U)$
are discrete, arithmetic subgroups of $\mathbb{H}^d/F^{*}$, and
$\mathbb{H}$ is the Hamiltonian algebra. We choose the $U$ carefully
so that the $\Gamma^i(U)$ are trivial. Let $\xi = p^{-
\sum v_i} T_p$, where the $v_i$ are scalars. Let $D(\mathbf{k},\alpha)$ be the number of eigenvalues of slope $\alpha$
of $\xi$ acting on $\mathcal{S}_{ \mathbf{k}}^D(U) $. Then, one of our main results is:

\begin{theorem}
Suppose $U$ is such that each $\Gamma^i(U)$ is trivial. There exist constants $\beta_1$ and $\beta_2$,
depending on $U$ such that if $\mathbf{k},\mathbf{k'} > n(\alpha) := [(\beta_1 \alpha - \beta_2)^d]$ and if $\mathbf{k} \equiv \mathbf{k'} \mod  p^{n(\alpha)}$,
then $D(\mathbf{k},\alpha) = D( \mathbf{k'}, \alpha)$
\end{theorem}

Let $K$ be an imaginary quadratic field of class number one and let
$O$ be the ring of integers of $K$. Fix an odd rational prime $p$,
which is inert in $K$.

We have an adelic definition of automorphic forms for modular forms
over number fields and in our case for imaginary quadratic fields
$K$. Let $K_{A}$ (resp $K_{A}^{*}$) be the adele ring (resp the
idele group) of $K$. We put $G_{K} = GL_{2}(K)$ and $G_{A} =
GL_{2}(K_{A})$. The center $Z_{A}$ of $G_A$ is isomorphic to
$K_{A}^{*}$. For a unitary character $\chi$ of the idele class group
$K_{A}^{*}/K^{*}$, let $L^{0}_{2}(G_{K} \setminus G_{A}, \chi)$
denote the space of measurable functions on $G_{A}$ satisfying
certain boundedness conditions.

We use the Eichler-Shimura relation between a subspace of
forms of this type and group cohomology. Essentially, if $\Gamma$ is
a congruence subgroup of $SL_2(O)$, and $S_{g,g}(O) = S_n(O) \otimes
S_g(O)$, where $S_g(O)$ is the $g$-th symmetric tensor power with an
action of $GL_2(O)$ on both components, then a certain space of
modular forms is isomorphic to $H^1(\Gamma, S_{g,g}(O))$. It is this
cohomology group that is useful to us. 

Let $m$ be the minimal number
of generators of $\Gamma$ and let $D(g,\alpha)$ denote the number of eigenvalues of slope $\alpha$ and weight $g$ for the $T_p$  operator. Our main theorem is:

\begin{theorem}
$D(g,\alpha)$ has an upper bound which is independent of the slope $\alpha$ and is
always less than $ [3m(\alpha+1)^2/2] m $.

\end{theorem}

\underline{Acknowledgements:}

The author would like to thank Fred Diamond for his supervision and
guidance for the last few years and Kevin Buzzard for a very careful
reading of the earlier drafts and for his suggestions. We also thank the referee for his/her comments.

\section{Preliminaries}

We refer the reader to the book \cite{DS} for a more detailed explanation of classical modular forms, cusp forms,
congruence subgroups and Hecke operators.

Let $M_k(SL_2 (\mathbb{Z}))$ denote the space of modular forms of weight $k$ for $SL_2 (\mathbb{Z})$, and for $\Gamma$ a congruence subgroup, let $M_k(\Gamma), S_k(\Gamma)$ denote the spaces of modular forms and cusp
forms of weight $k$ for $\Gamma$.

\begin{definition}

For congruence subgroups $\Gamma_1 , \Gamma_2$ of $SL_2(\mathbb{Z}),
 \alpha \in GL_2^+(\mathbb{Q})$, the weight $k$ Hecke operator $\Gamma_1 \alpha
 \Gamma_2$ operator takes functions $f \in M_k({\Gamma_1})$ to
 $M_k(\Gamma_2)$ (and similarly for cusp forms) by :

$\Gamma_1 \alpha \Gamma_2 : f \mapsto \sum_j f|_k \beta_j$ , where
$\beta_j$ are orbit representatives.

\end{definition}

Let $p$ be a prime. When $\alpha = \left(
                  \begin{array}{cc}
                    p & 0 \\
                    0 & 1 \\
                  \end{array}
                \right) $, we call the double coset operator the
                $T_p$ operator.

\begin{definition}

If $f$ is a non-zero modular form and an eigenform for $T_n$ for all
$n$, then it is called a Hecke eigenform. If $p$ is a prime, then
the $p$-adic valuation of the corresponding eigenvalue is called the
slope of an eigenform.

\end{definition}
We now define $p$-adic families of modular forms (note that they are not the same as $p$-adic modular forms).

\begin{definition}
Let $c \in \mathbb{Z}_p$ and for $r \geq 0$, let $B(c,r) = \{ k \in
\mathbb{Z}_p : |k - c| < r \}$. Let $N$ be an integer prime to $p$.
Then a $p$-adic family of modular forms of level $N$ is a formal
power series:

\center{$\sum_{n \geq 0} F_n q^n$},

where each $F_n : B(c,r) \rightarrow \mathbb{C}_p$ is a $p$-adic
analytic function, with the property that for all sufficiently large (rational) integers $k$, each $\sum F_n(k) q^n$ is the
Fourier expansion of a modular form of weight $k$.
\end{definition}

An example of a non-cupsidal family is the $p$-adic Eisenstein series $E_k^*(z) = E_k(z) -
p^{k-1}E_k(pz)$.

The slope $\alpha$ subspaces have been of great interest for a while,
and Gouvea and Mazur \cite{GM} made some very precise conjectures
about the dimensions of these spaces. Let $d(k,\alpha)$ be the
dimension of the slope $\alpha$ subspace of the space of classical
cuspidal eigenforms for the $T_p$ operators. Then the exact
conjecture (Buzzard and Calegari \cite{BC} found a counterexample a
few years ago) was:

\begin{itemize}
  \item If $k_1 , k_2 > 2\alpha + 2$,
  \item and $k_1 \equiv k_2 \mod p^n(p-1)$
  \item Then, $d(k_1, \alpha) = d(k_2, \alpha)$ (this condition is called local constancy).
\end{itemize}

Our goal is to get a Gouvea-Mazur type of result in the case
of Hilbert modular forms. In the case of modular forms over
imaginary quadratic fields we get an upper bound. To get these
results we need to define two more objects - Newton Polygons and Symmetric
tensor powers - as we use them both in the next two sections.

\smallskip

\underline{Newton Polygon:}

Let $L$ be a finite free $\mathbb{Z}_p$-module equipped with a
$\mathbb{Z}_p$-linear endomorphism $\xi$, so we can think of $L$ as a
$\mathbb{Z}_p [\xi]$ module. Let $\sum_{s=o}^{t} c_{s}X^{t-s}$ be the
characteristic polynomial of $\xi$ acting on $L \otimes \mathbb{Q}_{p}$.
Then, if $v_{p}$ denotes the usual $p$-adic valuation on $\mathbb{Z_{p}}$, we plot
the points $(i, v_{p}(c_{i}))$ in $R^{2}$, for $0 \leq i \leq t$,
ignoring the $i$ for which $c_{i}= 0$. Let $C$ denote the convex
hull of these point. The Newton polygon of $\xi$ on $L$ is the lower
faces of $C$, that is the union of the sides forming the lower of
the two routes from $(0,0)$ to $(t, v_{p}(c_{t}))$ on the boundary
of $C$. This graph gives us information about the $p$-adic
valuations of the eigenvalues of $\xi$. If the Newton polygon has a
side of slope $\alpha$ and whose projection onto the $x$ axis has
length $n$, then there are precisely $n$ eigenvectors of $\xi$ with
$p$-adic valuation equal to $\alpha$. The exact statement is:

\begin{theorem}

Let $\overline{L}$ be a field which is complete with respect to a valuation
$v$. Let $f(x) = \sum_{j=0}^{d}a_j x^j \in L[x]$ be a polynomial
with $a_0.a_d \neq 0$. Let $l$ be a line segment of the Newton
polygon of $f$ joining $(j,v(a_j))$ and $(h,v(a_h))$ with $j < h$.
Then $f(x)$ has exactly $h-j$ roots $\gamma$ in $L$ such that
$v(\gamma)$ is the negative of the slope of $l$.

\end{theorem}

\underline{Symmetric tensor powers:}

Let $R$ be any commutative ring. For any $R$-algebra $A$ and for
$a,b \in \mathbb{Z}_{\geq 0}$, we let $S_{a,b}(A)$ denote the $M_2(R)$-module
$Symm^a(A^2)$ (the $a^{th}$ symmetric power with $M_2(R)$ action). The action is
given by $x\alpha = (det \alpha)^{b}xS^{a}(\alpha)$. If $A^2$ has a natural basis $e_1, e_2$, then $S_{a,b}(A)$ has a
basis $f_0,...,f_a$ where each $f_i = e_1^{\otimes i} \otimes
e_2^{\otimes (a-i)}$

Another way to think of this action is to use the equivalence of
$Symm^{a}(A^2)$ with the space of homogeneous polynomials of degree
$a$ in $2$ variables. The $M_{2}$ action can be described as follows.

Let $A$ be an $R$-algebra and consider the polynomial ring $A[x,y]$.
Let $f(x,y) \in A[x,y]$. If  $\alpha = \bp a &b\\c&d\ep$, then
$\alpha f(x,y) = f(ax+cy, bx+dy)$. We can also twist this
action by the determinant, $\alpha f(x,y) =
det(\alpha)^{b}f(ax+cy, bx+dy)$. When considered as a module for $SL_2(\mathbb{R})$ we know $det =
1$, so we'll call the $g$-th symmetric tensor power $S_g$.

$\mathbf{Strategy}$

Let $L$ be a finite, free $R$-module of rank $t$ ($L$ will correspond to a space of automorphic forms). 
We will define $K$ to be a submodule of $L$
such that $L/K \equiv \oplus O/p^{a_{i}}O$ with the $a_{i}$
decreasing and $a_{i} \leq n$. We consider the characteristic polynomial $p(x)$ of $\xi$ acting on $L$ and plot
its Newton polygon. 

We let $L'$ be a space of forms corresponding to a different weight, and choose $K'$ to be a submodule similar to $K$.
 Modulo a certain power of $p$ we show that the
spaces $L/K$ and $L'/K'$are isomorphic which leads to congruences of the
coefficients of the respective characteristic polynomial. 
This tells us that the Newton polygons of fixed slopes coincide which gives us
local constancy of the slope $\alpha$ spaces.

In the case that the quotient spaces are not isomorphic, we get a lower bound for
the Newton polygon associated to $\xi$ which transforms into an
upper bound for the slope $\alpha$ eigenspaces.

\section{The Totally real case}

\subsection{Automorphic forms over Quaternion Algebras}

In this section we'll define automorphic forms over quaternion
algebras. Due to Jacquet Langlands, they correspond to a subspace of cuspidal Hilbert modular forms. 
The advantage is that one can work
more easily with the definitions for quaternion algebras as things
are finite. For a more detailed description see \cite{T2} and \cite{Hi1}.

If $F$ is a field, then a quaternion algebra $D$ over $F$ is a
central, simple algebra of dimension $4$ over $F$. Central means that
$F$ is the center of $D$ and simple means that there are no
two-sided ideals of $D$ except for $\{0\}$ and $D$ itself. For each
embedding $\sigma : F \hookrightarrow \mathbb{R}$, we say that $D$ is
ramified at $\sigma$ if $D_{\sigma} = D \otimes_{F, \sigma}
\mathbb{R} \cong \mathbb{H}$,
 where $\mathbb{H}$ is the Hamilton quaternion algebra. A totally
definite quaternion algebra means that it ramifies at exactly all the
infinite places. Let $\mathbb{A}$ be the ring of adeles.

Let $F$ be a totally real field, $[F:\mathbb{Q}] = d$, where $d$ is
even. Let $K$ be a Galois extension of $\mathbb{Q}$, which splits $D$,
with $F \subseteq K$. Fix an isomorphism $D \otimes_F K \cong M_2
(K)$. Assume that $D$ is a totally definite quaternion algebra over $F$,
unramified at all finite places and fix $O_{D}$ to be a maximal
order of $D$.

Let $G = Res_{F/\mathbb{Q}} D^{*}$ be the algebraic group defined by
restriction of scalars. Fix $\mathbf{k} = (k_{\tau}) \in \mathbb{Z}^{I}$ such that each component $k_{\tau}$
is $\geq 2$ and all components have the same parity. Set
$\mathbf{t} = (1,1,...,1) \in Z^{I}$ and set $\mathbf{m} = \mathbf{k} - 2\mathbf{t}$. 
Also choose $\mathbf{v} \in
Z^{I}$ such that each $v_{\tau} \geq 0$, some $v_{\tau} = 0$ and $\mathbf{m}
+ 2\mathbf{v} = \mu \mathbf{t}$ for some $\mu \in \mathbb{Z}_{\geq 0}$.

For any $R$-algebra $A$ and for $a,b \in \mathbb{Z}_{\geq 0}$, we let
$S_{a,b}(A)$ denote the left $M_2(R)$-module $Symm^a(A^2)$ (with the $M_2(R)$ action described in the previous section). 
If $\mathbf{k} \in \mathbb{Z}[I]$ and $\mathbf{m}, 
\mathbf{v}, \mu$ are as before we set 
$L_{\mathbf{k}} = \otimes_{\tau \in I} S_{m_\tau, v_\tau }(\mathbb{C})$. If $R$ is a ring such that $O_{K,v} \subseteq R$, for some $v | p$,
then,  $L_{\mathbf{k}}(R) = \otimes_{\tau \in I} S_{m_{\tau},v_{\tau}}(R) $

Now, we'll define automorphic forms on these quaternion algebras.

First, we pick a prime rational $p$ which is inert in $K$. Let $M$ be the semigroup in $M_2(O_{F,p})$ consisting of matrices
$\bp a &b\\c&d\ep$ such that $c \equiv 0 \mod p$ and $d \equiv 1
\mod p$.  Let $U \subseteq G_{f}$ be
an open compact subgroup such that the projection to $G(F_p)$ lies
inside $M$. If $u \in U$, let $u_{p} \in G(F_{p})$  denote the image under the
projection map.

Next, we define a weight $\mathbf{k}$ operator.

If, $ f: G(\mathbb{A}) \rightarrow L_{\mathbf{k}}(R)$ and $u =
u_{f}.u_{\infty} \in G(\mathbb{A})$ then,

$(f|_{\mathbf{k}}u)(x) = u_{\infty}f(x.u^{-1})$, when $R = \mathbb{C}$.

$(f||_{\mathbf{k}}u)(x) = u_{p}f(x.u^{-1})$, when $R$ is an $O_{K,p}$-algebra.

The space of automorphic forms for $D$, of level $U$ and weight $\mathbf{k}$ can be described as:

$\mathcal{S}_{\mathbf{k}}^D(U) = \{ f: D^{*} \setminus G(\mathbb{A})
\rightarrow L_{\mathbf{k}} \mid f |_{\mathbf{k}}u = f , \forall u \in U\}$ $= \{f: G_{f}/U \rightarrow L_{\mathbf{k}} \mid f(\alpha.x) = \alpha.f(x),
\forall  \alpha \in D^{*} \}$

\smallskip

$\mathcal{S}^{D}_{\mathbf{k}}(U, R) = \{ f: D^{*} \setminus G(\mathbb{A})
\rightarrow L_{k}(R) \mid  f ||_{\mathbf{k}} u = f , \forall u \in U \}$.

\smallskip

The purpose of introducing $\mathcal{S}^{D}_{\mathbf{k}}(U, R)$ is to give
$\mathcal{S}^{D}_{\mathbf{k}}(U)$ an integral structure which allows us to think of $\mathcal{S}^{D}_{\mathbf{k}}(U, R)$ as $\oplus_{\gamma_i \in X(U)} (\gamma_i  L_k(R))^{D^{*} \cap \gamma_i
U \gamma_i^{-1}}$. Thus, we see that $S_{\mathbf{k}}^{D}(U,R)$ is an $R$-lattice in
$S_{\mathbf{k}}^{D}(U)$.

\smallskip

Let $X(U) = D^{*}\backslash G_{f} / U$. We know this is finite, so let $h = |X(U)|$
and let $\{\gamma_{i} \}_{i =1}^{h}$ be the coset
representatives. So, $G_{f} = \coprod_{i = 1}^{h} D^{*}.\gamma_{i}.U$.

Define $\widetilde{\Gamma^{i}(U)} := D^{*} \cap
\gamma_{i}.U.G_{\infty}^{D}.\gamma_{i}^{-1}$ and let $\Gamma^{i}(U) := \widetilde{\Gamma^{i}(U)}/\widetilde{\Gamma^{i}(U)} \cap F^{*}$.

We want to impose conditions on $U$ such that the
${\Gamma^{i}(U)}$ are trivial. We know that $\widetilde{\Gamma^{i}(U)}$
are discrete arithmetic subgroups of $G_{\infty}$ and that
${\Gamma^{i}(U)}$ are discrete in
$G_{\infty,+}/F_{\infty}^{*}$. As $D$ is a totally definite quaternion algebra  $G_{\infty} \cong
(\mathbb{H})^{d}$, where $\mathbb{H}$ is the Hamiltonian algebra. So, ${\Gamma^{i}(U)}$ is discrete in
$G_{\infty,+}/F_{\infty}^{*}$ which is compact. Thus,
${\Gamma^{i}(U)}$ is finite. 

Let $N$ be an ideal in $O_{F}$. Define:

$U_{0}(N) = \{ \left(
                 \begin{array}{cc}
                   a & b \\
                   c & d \\
                 \end{array}
               \right)
 \in \prod_{q} GL_{2}(O_{F,q}) | c_q \in N O_{F,q} \forall q \}$, where $q$ runs over all the
finite primes of $F$.

$U_{1}(N) = \{ \bp a &b\\c&d\ep \in U_{0}(N) \mid  a - 1 \in N
\}$,

We cite a result by Hida \cite{Hi1} (Sec $7$).

\begin{lemma}(Hida)

Put $U(N) = \{ x \in U_1(N) : x_N = \bp a &b\\c&d\ep , d - 1 \in N
\}$ for
 each ideal $N$ of $O_F$ . Let $l$ be a prime ideal of $O_F$ and let $e$ be the
  ramification index of $l$ over $Q$. Then, if $s > 2e/(l-1)$, then ${\Gamma^{i}(U(l^s))}$
   is torsion free for all $i$

\end{lemma}

According to this result $U$ can be chosen such that
${\Gamma^{i}(U)}$ are torsion free for all $i$. Coupled
with the statements above, this means that
${\Gamma^{i}(U)}$ are trivial, provided $U$ is chosen
carefully.

\smallskip

Hecke operators are defined in a similar fashion as in classical modular forms using the double coset
decomposition. Let $U,U'$ be open compact subgroups and $x \in G_{f}$. We define:

$\xi = UxU': \mathcal{S}^{D}_{\mathbf{k}}(U,R) \rightarrow
\mathcal{S}^{D}_{\mathbf{k}}(U',R)$, where

$\xi : f \mapsto \sum f ||_{\mathbf{k}} x_i$, where $UxU = \coprod U x_{i} $.

In particular, we have the $T_q$ operator. $T_{q} = [U \eta_{q} U] $, where $\eta_{q} = \bp \pi_{q} &0\\0&1\ep
$ and $\pi_{q} \in F_{f}$ is $1$ everywhere except at $q$, where it
is a uniformiser.

\subsection{Calculations}

\smallskip

Let $n \leq k_{i}$.

We can think of $L_{\mathbf{k}}(R) = L_{k_{1}}(R) \otimes L_{k_{2}}(R)
\otimes ...... \otimes L_{k_{d}}(R)\otimes det()^{n_k}$, where $det()^{n_k}$ accounts for the twist by determinants and
$L_{k_i}(R)$ are simply the $k_i$th symmetric powers.

 We define $W_{\mathbf{k}}^{n}(R)$ to be generated by the submodules
$W_{k_{1}}^{n}(R) \otimes L_{k_{2}}(R) \otimes ...... \otimes
L_{k_{d}}(R)\otimes det()^{n_k}$ , $L_{k_{1}}(R) \otimes
W_{k_{2}}^{n}(R) \otimes ...... \otimes L_{k_{d}}(R)\otimes
det()^{n_k}$ and up to $L_{k_{1}}(R) \otimes L_{k_{2}}(R) \otimes ...... \otimes
W_{k_{d}}(R)\otimes det()^{n_k}$, where each $W_{k_{i}}^{n}(R)$ is generated by the $(n+1)$ $R$
submodules $\{ p^{n-j} x^{j} L_{k-j}(R) \}_{j=0}^{n}$.

We can think of each $W_{k_{i}}^{n}(R)$ in another manner. Note that $R[x,y] =
\oplus L_{k_{i}}(R)$. Let $J \subseteq R[x,y]$ denote the
homogeneous ideal $(p,x)$. Then for all $n > 0$, $J^{n}$ is also
homogeneous, so it can be written as $J^{n} = \oplus
W^{n}_{k_{i}}(R)$. It is not hard to check that the
$W_{k_{i}}^{n}(R)$ are invariant under $\Gamma_{j}(U)$. The key is
that if $\bp a &b\\c&d\ep \in \Gamma^j (U)$ then $c \equiv 0 \mod
p$.

\smallskip

\underline{Case 1: ${\Gamma^{i}(U)}$ are trivial:}

\smallskip

Thus, $\mathcal{S}^{D}_{\mathbf{k}}(U,R)/ \mathcal{W}_{\mathbf{k}}^{D}(U,R) \cong
\oplus_{i=1}^{h}L_{\overrightarrow{k}}(R)/\oplus_{i=1}^{h}W^{n}_{\mathbf{k}}(R)$ for all
$\mathbf{k}=(k_{1},k_{2},...,k_{d})$

Now, we use the fact that $L_1 \otimes L_2 /C \cong L_1 / W_1
\otimes L_2 / W_2 $ ,where $C = W_1 \otimes L_2 + L_1 \otimes W_2$. 

Let $I_{k_{i}^{n}}(R) :=  L_{k_i} / W_{k_i}^{n} \cong
\oplus_{j=1}^{n}O/p^{j}O$. This quotient depends on $n$.
 
Thus, $L_{\mathbf{k}}(R)/W_{\mathbf{k}}^{n}(R) \cong I_{k_1}^n \otimes ... \otimes
I_{k_d}^n  = I^{n}_{\mathbf{k}}(R) \cong \underbrace{\oplus (O/p^{n}O)}_{\sigma_1 \
\hbox{\scriptsize times}} \underbrace {\oplus
(O/p^{n-1}O)}_{\sigma_2 \ \hbox{\scriptsize times}} .....
\underbrace{\oplus (O/p^{1}O))}_{\sigma_{n} \ \hbox{\scriptsize
times}}$, where $\sigma_{i}$ is the multiplicity of each factor.

$\cong \oplus_{i=1}^r O/p^{a_{i}}O$, and where $a_1 \geq a_2 \geq
...$.

\smallskip

\underline{Case 2 ${\Gamma^{i}(U)}$ are not trivial:}

In this case, we see $\mathcal{S}_{\mathbf{k}}^{D}(U,R) =
\oplus_{i=1}^{h}(L_{\mathbf{k}}(R))^{\overline{\Gamma^{i}(U)}}$.

Let $\mathcal{W}_{\mathbf{k}}^{D}(U,R) = \{ f \in
\mathcal{S^{D}}_{\mathbf{k}}(U,R) | f(x) \in W_{\mathbf{k}}^{n}(R) \} \cong
\oplus_{i=1}^{h}W_{\mathbf{k}}^{n}(R)^{{\Gamma^{i}(U)}}$.

Thus, $\mathcal{S}^{D}_{\mathbf{k}}(U,R)/ \mathcal{W}_{\mathbf{k}}^{D}(U,R)
\cong \oplus_{i=1}^{r} O/p^{a_{i}^{k}}O$ for each
$\mathbf{k}=(k_{1},k_{2},...,k_{d})$.

\smallskip

We define functions $B, T$ based on the $\sigma_{i}$. These functions serve as upper and lower bounds for the
Newton polygons of the Hecke operators that are central to the results.

Let $b_{i} = n - a_{i}$ and $B(j) = \sum_{i=1}^j b_i$, so we get the following formulae 

\[ \ba{rcl}   \\
  b_{i}& = &\left\{ \ba{ll} j , & \mbox{if } (\sum_{k=1}^{j}\sigma_{k})h \leq i \leq
(\sum_{k=1}^{j+1}\sigma_{i})h$ , $0 \leq j \leq n \\
n, &\mbox{if } i \geq (\sum_{k=1}^{n}\sigma_{k})h . \ea \right. \ea
\]

\[ \ba{rcl} B(x): \mathbb{R} &\to& \mathbb{R}  \\
  B(x)& = &\left\{ \ba{ll}  \sum_{k=1}^{j}\sigma_{k} (k-1) + j(x -
\sum_{k=1}^{j}\sigma_{k}) , & \mbox{if } (\sum_{k=1}^{j}\sigma_{k})h
\leq x \leq
(\sum_{k=1}^{j+1}\sigma_{k})h , 0 \leq j \leq n \\
{\sum_{k=1}^{n}\sigma_{k}(n-1)} + n(x - \sum_{k=1}^{n}\sigma_{k}),
&\mbox{if }x \geq (\sum_{k=1}^{n}\sigma_{k})h. \ea \right. \ea \]

Now, let $M$ be the smallest integer such that $2M \geq n$, let
$T(x) = M + B(x-1)$. We see that $T$ can be described as:

\[ \ba{rcl} T(x+1): \mathbb{R} &\to& \mathbb{R}  \\
  T(x+1)& = &\left\{ \ba{ll} M , & \mbox{if } i = 0 \\ M + \sum_{k=1}^{j}\sigma_{k}
 (k-1) + j(x -
\sum_{k=1}^{j}\sigma_{k}) , & \mbox{if } (\sum_{k=1}^{j}\sigma_{k})h
\leq x \leq
(\sum_{k=1}^{j + 1}\sigma_{k})h , 0 \leq j \leq n \\
M + {\sum_{k=1}^{n}\sigma_{k}(n-1)} + n(x -
\sum_{k=1}^{n}\sigma_{k}), &\mbox{if }x \geq
(\sum_{k=1}^{n}\sigma_{k})h. \ea \right. \ea \]

Let $c = \inf \{T(x)/x\}$ for $x \geq 1$.

We consider the case when the $\overline{\Gamma^i}$ are
trivial, because it allows us to get a precise value for each of the $\sigma_i$
 in the above formula.

Since the structure of  $L_{\overrightarrow{k}}(R)/W_{\overrightarrow{k}}^{n}(R)$ is the tensor product of $d$ copies of
$\oplus_{i=1}^n O/p^i O$, we can see that $\sigma_1 = 1.h$ as we get
the $n$-th power only once. $\sigma_2 = (2^d-1)h$. Carrying on we
see that $\sigma_i = (i^d - (i-1)^d)h$ which gives us that $\sum_{i=1}^j \sigma_i = j^d h $. So, $B(x)$ is a piecewise linear function which has slope
 $r$ for $ r^d h \leq x \leq (r+1)^d h$.

Now, $\inf T(x)/x = \inf_{x \ge 0} P(x)/(x+1)$ which is at least $P(x)/2x$
for $x \ge 1$ and $P(0)/2$ for $x \le 1$.

Let $q(x) = (x/h)^{1/d} - 1 $. It follows that $q(x) < B'(x)$. If
$Q(x) = \int_0^y q(y)dy$, then $Q(x) < B(x)$.

$Q(x) = (\frac{d}{d+1})(\frac{x}{h})^{(d + 1)/d} - x$.  Let $P(x) = M + Q(x) < T(x + 1)$. So, $\frac{P(x)}{x} = \frac{M}{x} + \frac{d}{d+1}
(\frac{x}{h})^{1/d} - 1$.

To find the minimum of $P(x)/x$, we use basic calculus.

$(\frac{P(x)}{x})'= 0$ $\Rightarrow \frac{-M}{x^2} + \frac{1}{d+1} (\frac{x}{h})^{1/d -1} =
0$ $\Rightarrow (\frac{x}{h})^{d+1/d} = (d+1)M$. Therefore, the minimum of $P(x)/x$ occurs at $x = h(M(d+1))^{d/(d+1)}$.
Putting this value in for $P(x)/x$, we get that

$\frac{P(x)}{x} \geq \frac{M^{1/(d+1)}}{(h(d+1))^{d/(d+1)}} +
\frac{d}{(d+1)} ((d+1)M)^{1/(d+1)} - 1$. So, $\frac{P(x)}{x} \geq M^{1/(d+1)}(\frac{1}{d+1})^{d/(d+1)}$

Recall that we chose $M$ such that $2M \geq n$ which means that $\frac{P(x)}{x} > c_1 n^{1/(d+1)} - 1$,
 where $c_1 = (\frac{1}{d+1})^{d/(d+1)}(\frac{1}{h^{d/(d+1)}} + 1)$. Note that this is
 true only for $x \leq n^d (h + 1)$. Since $c = \inf \{T(x)/x\}$, we see that
 $c = min \{ c_1 n^{1/(d+1)}, n \}$.

\smallskip

In the preliminaries  section we outlined our strategy, where the goal was to show that the Newton polygons of certain spaces coincide. We will now define our $L$, $K$ and show the exact proposition (due to Buzzard) that we use later.

Let $L$ be a finite, free $R$-module of rank $t$ (where $L$ corresponds to $\mathcal{S}_{\overrightarrow{k}}^{D}(U,R)$ ),
equipped with an $R$-linear endomorphism $\xi$.
 Define $K$ to be a submodule ($K$ is $\mathcal{W}_{\overrightarrow{k}}^{D}(U,R)$)
such that $L/K \cong \oplus O/p^{a_{i}}O$ with the $a_{i}$
decreasing and $a_{i} \leq n$. Let $b_{i} = n - a_{i}$. Let $M$ be
the smallest integer such that $2M \geq n$. Let $p(x) = \sum_{i=1}^{t}d_{i}x^{i}$ be the characteristic
polynomial of $\xi$ acting on $L$.

Then, we have the following result.

\begin{proposition}
\begin{enumerate}
    \item If the above conditions hold, then the Newton polygon of $\xi$ lies above the
    function $B$.

\item If $\alpha < c$, and $K \subset L $ and $K' \subset L$ satisfy the above hypothesis, and if $L/K \cong L'/K'$
as $R[\xi]$-modules, then the Newton Polygons of $\xi$ of small slope will coincide.
\end{enumerate}
\end{proposition}

\begin{proof}

\begin{enumerate}
    \item Choose an $R$-basis $(e_{i})$ for $L$ such that
$(p^{a_{i}}e_{i})$ is an $R$-basis for $K$. Let $(u_{i,j})$ be the
matrix of $\xi$ acting on $L$. As $\xi(K) \subseteq p^{n}(L)$ we get
that $p^{b_{j}}$ divides $u_{i,j}$.

By the definition of the characteristic polynomial $det (x -
u_{i,j}) = \sum d_s x^{t - s}$. Then, we can see that:

$d_s = (-1)^{s} \sum_{J \subseteq \{ 1,2,...,t\} \text{of size s}}
\sum_{\sigma \in Symm(J)} sgn(\sigma) \prod_{j \in J} u_{j, \sigma
(j)}$.

Now, we know that $p^{b_{\sigma j}} $ divides $u_{j, \sigma j}$, and
so we get that:

$d_s = (-1)^{s} \sum_{J \subseteq \{ 1,2,...,t\} \text {of size s}}
p^{\sum_{j \in J} b_j} \sum_{\sigma \in Symm(J)} sgn(\sigma)
\prod_{j \in J} u_{j, \sigma j}/p^{b_{\sigma(j)}}$.

Since, the $b_i$ are increasing and for all $J$ of size $s$,
$\sum_{j \in J} b_j \geq b_1 + b_2 + ... + b_s $. This gives us $d_s$ is divisible by $p^{B(s)}$ for all $s$. This
means that $d_{i} \equiv 0 \mod p^{B(i)}$ which gives us the first
statement.

 \item Let $\sum d_{i}x^{i}$ be the characteristic polynomial of $\xi$ on $L$.
Let $L'$ be another free $R$ module of rank $t'$, and $K'$ a
submodule such that $L'/K' \cong L/K$ and $\xi(K') \subseteq
p^{n}L'$. Let $\sum d'_{t -i}x^{i}$ be the characteristic polynomial
of $\xi$ on $L'$. Set $d'_i = 0$ for $i > t$. Assume $t' < t$.

Claim: $d_{i} \equiv d'_{i} \mod p^{T(i)}$.

Proof of claim:

Since $\sum d_{t - s} x^s = det (x.I - (u_{j,k}))$, we can expand
the coefficients
 $d_i$ in terms of the matrix coefficients $u_{i,j}$ as follows.

$d_s = (-1)^s \sum_{J \subseteq \{1,2,...t\} } \sum_{\sigma \in
symm(J)} sgn(\sigma)
 \prod_{j \in J}  u_{j, \sigma(j)}$, where $J$ is a set of size $s$.

Let $(e_{i})$ be a basis of $L$ such that $p^{a_i}e_i$ is a basis of
$K$. Let $f_i$ be the reduction of each $e_i$ in $L/K =
\overline{L}$. Since $\overline{L} \equiv \overline{L'}$, choose
$f'_i$ in $\overline{L'}$ through this isomorphism and let $e'_i$ be
the lift to $L$ of each of these $f'_i$.
 Then $p^{a_i}e'_i$
is a basis for $K'$. Let $(u_{j,k})$ be the matrix for $\xi$ acting on $L'$. Since we
know that $L'/K' \cong L/K$ we can infer that $u_{j,k} \equiv
u'_{j,k} \mod p^{a_j}$. We set $u'_{j,k} = 0$, if $max \{j,k \} >
t'$. To establish the claim we need to show that $\prod u_{j,k} \equiv
\prod u'_{j,k} \mod p^{T(s)}$. We'll show $\prod_{k} u_{j,k} \equiv \prod_{k} u'_{j,k} \mod
p^{N_J}$, where $N_J \geq T(s)$.  

$u_{j,k} \equiv u'_{j,k} \mod p^{a_j}$ $\Rightarrow u_{j,k}/p^{b_k} \equiv u'_{j,k}/p^{b_k} \mod
p^{c_{j,k}}$, where $c_{j,k} = max \{ a_j - b_k, 0 \}$

$\Rightarrow \prod_{j \in J} u_{j,\sigma(j)}/p^{b_{\sigma(j)}}
 \equiv \prod_{j \in J} u'_{j,\sigma(j)}/p^{b_{\sigma(j)}} \mod p^{ min_{j \in J} \{ c_{j,\sigma(j)} \} }$ $\Rightarrow \prod_{j \in J} u_{j,\sigma(j)} \equiv \prod_{j \in J}
u'_{j,\sigma(j)} \ \mod p^{N_J}$, where $N_J = \sum_{j \in J} b_j + min_{j \in J}\{c_{j,\sigma(j)}\}$ $\geq \sum_{j \in J} b_j + c_{j_{0},\sigma (j_0)}$, where $j_0 = max_{j \in
J} j$ $\geq \sum_{j_0 \neq j \in J} b_j + b_{j_0} + c_{j_{0},\sigma (j_0)}$ $\geq \sum_{j_0 \neq j \in J} b_j + max \{a_{j_0}, b_{j_0}\}$ $\geq b_1 + b_2 + .... + b_{s-1} + M $ $ = T(s)$.

Now, the function $T(i)$ is convex, piecewise linear and $c <
T(i)/i$. So, $\alpha < c \Rightarrow \alpha i < c i < T(i)$. This says that if the Newton Polygon has a side of slope $\alpha <  c$, then it
lies below the graph of $T$. The endpoints of this side are $(s_1,
v_{p}(T(s_1)))$ and $(s_2,(T(s_2)))$. So, we see that
$v_p(d_{s_i}) < T(s_i)$ for $ i = 1,2$. Now, $d_{i} \equiv d'_{i} \mod p^{T(i)}$ we conclude
 that a Newton Polygon of this length and slope
depends only on $\overline{L}$. So, we get that for $\alpha < c$ the number of eigenvalues of $\xi$
on $L$ with slope $\alpha$ depends only on the isomorphism class of
L.

\end{enumerate}

\end{proof}

\subsection{Lemmas}

Fix a prime $p$ in $\mathbb{Q}$, which is inert in $K$. Let $R =
O_{K,p}$.

We know that $\mathcal{S}^{D}_{\mathbf{k}}(U,R) =
\oplus_{i=1}^{h}(L_{\mathbf{k}}(R))^{{\Gamma^{i}(U)}}$.

Now, we have two cases.

\begin{itemize}
    \item ${\Gamma^{i}}(U)$ are trivial.
\item ${\Gamma^{i}}(U)$ are not trivial.
\end{itemize}

\smallskip

\begin{lemma}
Let $\xi = U\eta_{p}U$, $U=U_{1}(Np)$, then
$\xi(\mathcal{W}_{\mathbf{k}}^{D}(U,R)) \subseteq p^{n + \sum v_i}
\mathcal{S}^{D}_{\mathbf{k}}(U,R) $, where $\eta_p = \bp p &0\\0&1\ep $, and the $v_i$ are constants chosen in the definition of $\mathcal{S}^{D}_{\mathbf{k}}(U,R) $
\end{lemma}

\begin{proof}
First, note that if $\xi = U\eta_{p}U = \coprod U \eta_{i} $, then
$(\eta_{i})_{p}\equiv \bp 0 &*\\0&*\ep \mod p$.

If $f \in \mathcal{W}_{\mathbf{k}}^{D}(U,R)$, then $\xi : f \mapsto \sum
 f ||_{\mathbf{k}} \eta_{i}$. Now, $(f ||_{\mathbf{k}}\eta_{i})(x) = (\eta_{i})_{p}f(x.\eta_{i}^{-1}) $,
where $f(x.\eta_{i}^{-1}) \in \mathcal{W}_{k}^{n}(U,R)$ i.e., $f(x.\eta_{i}^{-1}) = f_{1}(x_{1},y_{1})\otimes .....\otimes
f_{d}(x_{d},y_{d})$, where $f_{i}(x_{i},y_{i}) \in W_{k_{i}}^{n}(R)$
for some $i$.

Under the action of a matrix $\bp a &b\\c&d\ep$, where $a,c \equiv 0
\mod p$, we see $\bp a &b\\c&d\ep f_{i}(x_{i},y_{i}) = f_{i}(ax_{i} + cy_{i}, bx_{i}
+ dy_{i})$ is divisible by $p^{n}$, because if $\gamma = \bp a &b\\c&d\ep $ as above,
 then $\gamma x_i$ is divisible by $p$. So, $\gamma <p, x_i> \subseteq p R[x_i,y_i]$ which means that
$\gamma <p, x_i>^{n} \subseteq p^{n} R[x_i,y_i]$.

We also need to account for the twists by the determinant factor. As
we had written earlier, the action of $\gamma = \bp a &b\\c&d\ep $
on $f(x,y)$ is given by $det(\gamma)^{v_i} f(ax + cy, bx + dy)$,
where the $v_i$ were chosen while defining $\mathcal{S}^{D}_{\mathbf{k}}(U,R)$. Since $\gamma = \bp a &b\\c&d\ep$, where $a, c \equiv o \mod p$,
 the determinant factor gives us an extra power of $p$. Now, in all the other
  factors we get $p^{v_j}$. Since it is a tensor product of all these
  factors we see that $\xi (f)$ is divisible by $p^{n + \sum v_i}$. Thus, $\xi(\mathcal{W}_{\mathbf{k}}^{D}(U,R)) \subseteq p^{n + \sum
v_i}\mathcal{S}^{D}_{\mathbf{k}}(U,R) $.

\end{proof}

\begin{lemma}
Assume that ${\Gamma^{i}}(U)$ are trivial. If $\mathbf{k} \equiv \mathbf{k'} \mod
p^{n-1}$, then $\mathcal{S}^{D}_{\mathbf{k}}(U,R)/ \mathcal{W}_{\mathbf{k}}^{D}(U,R)
\cong \mathcal{S}^{D}_{\mathbf{k'}}(U,R)/ \mathcal{W}_{\mathbf{k'}}^{D}(U,R)$
\end{lemma}

\begin{proof}

We define $\mathbf{k'}$ as follows.

$\mathbf{k'} = \mathbf{k} + (0,0,....,0,p^{n-1},0,...,0)$. i.e.

\[ \ba{rcl}   \\
  k'_{i}& = &\left\{ \ba{ll} k_{i} , & \mbox{if } i \neq i_{0} \\
k_{i} + p^{n-1}, &\mbox{if } i = i_{0} . \ea \right. \ea
\]

Now, for each $k_{i} \equiv k_{i}' \mod p^{n-1}$, we want to show
that $I^{n}_{k_{i}}(R) \cong I^{n}_{k_{i}'}(R)$.

As abelian groups, each $I^{n}_{k_{i}}(R) \cong
\oplus_{j=1}^{n}O/p^{j}O$, which means that $I^{n}_{\mathbf{k}}(R) \cong
\oplus_{j=1}^{r}O/p^{a_{j}}O$.

We want to prove that as $M$ modules $I^{n}_{\mathbf{k}}(R) \cong I^{n}_{\mathbf{k'}}(R)$
 if $\mathbf{k} \equiv \mathbf{k'} \mod p^{n-1}$. Let $\phi : L_{\mathbf{k}}(R) \rightarrow L_{\mathbf{k}}(R)$, where $\phi$ is the identity on each component except $i_{0}$.

On $L_{k_{i_0}}$, we can think of it as $\phi : f(x,y) \mapsto
f(x,y)y_{i_0}^{p^{n-1}}$. (Note that $x = (x_1, x_2,....,x_d)$ and $y = (y_1,y_2,.....,y_d)$). To show that $\phi$ induces a homomorphism we need to verify that

$\phi(\gamma f(x,y)) = \gamma \phi(f(x,y))$, where $\gamma = \bp a
&b\\c&d\ep $ $\Rightarrow \phi(\gamma f(x,y)) - \gamma \phi(f(x,y)) \in W_{\mathbf{k} +
p^{n-1}}^{n}(R)$. 

We'll show it for each $i$, i.e., $\phi(\gamma_i f(x_i,y_i)) - \gamma_i
 \phi(f(x_i,y_i)) \in
<p,x_i>^n$

Recall that for each $i$, $<p,x_i>^n = \oplus W_{k_{i}}^{n}$. Now, $\phi(\gamma_i f(x_i,y_i)) - \gamma_i
 \phi(f(x_i,y_i)) $.

$= \phi(f(a_{i}x_{i} + c_{i}y_{i}, b_{i}x_{i} + d_{i}y_{i})) - \bp
a_i &b_i\\c_i&d_i\ep \phi(f(x_i,y_i))$
 $= (f(a_i x_i + c_i y_i, b_i x_i + d_i y_i))y_i^{p^{n-1}} - \bp a_i
&b_i\\c_i&d_i\ep (f(x_i,y_i))y_i^{p^{n-1}}$ 
$= (f(a_i x_i + c_i y_i, b_i x_i + d_i y_i))y_i^{p^{n-1}} - (f(a_i
x_i + c_i y_i,b_i x_i + d_i y_i))(b_i x_i + d_i y_i)^{p^{n-1}}$

$= (f(a_i x_i + c_i y_i, b_i x_i + d_i y_i))(y_i^{p^{n-1}} - (b_i
x_i + d_i y_i)^{p^{n-1}})$. It will be enough to show that $(y_i^{p^{n-1}} - (b_i x + d_i
y_i)^{p^{n-1}}) \in <p,x_i>^{n}$.

Claim:

If $f,f' \in L_{k_i}(R)$, such that $f - f' \in W_{k_i}^{n}$,
 then $f^{p^{s}} - f'^{p^{s}} \in W_{p^{s}k_i}^{n}$.

Proof of claim:

Let $f = f' + F$ for some polynomial $F \in <p,x_i>^n$. Then $f^p =
(f' + F)^p $ which means that $f^p = f'^p + G$, where $G \in
(pF,F^p) \subseteq <p,x_i>^n$. The rest of the claim follows by
induction.

\smallskip

Now, we show that $(y_{i}^{p^{n-1}} - (b_i x_i + d_i y_i)^{p^{n-1}}) \in
<p,x_i>^{n}$. For $n=1$, we have to know that $y_i - (b_i x_i + d_i y_i) =
((1-d_i)y_i - b_i x_i)$. Now, $d_i \equiv 1 \mod p$, so we're done. Then for $n > 1$, we use the claim with $n=1$ and the claim with $s
= n -1, n = 1$. Let $f = y_i, f' = (b_i x_i + d_i y_i)$. By the
previous step, we know that $f - f' \in <p,x_i> \subseteq
W_{k_i}^n$. So, by the claim above we know that $f^{p^{s}} -
f'^{p^{s}} \in W_{p^{s}k_i}^{n}$ i.e. $(y^{p^{n-1}} - (b_i x_i + d_i
y_i)^{p^{n-1}}) \in <p,x_i>^{n}$.

\end{proof}

$\mathbf{Summary}:$

We know that the Newton Polygon is described in terms of
$v_{p}(d_{i})$. The proposition says that it lies above the $B(i)$
and that line segments of the polygon lying below $T(i)$ depend only
on $L/K$. So, the functions $B,T$ form the bottom and top boundaries
of a region where the Newton Polygon depends only on the isomorphism
class of $L/K$. Therefore, if the eigenvalues of $\xi$ have small
slope then the side of the Newton polygon which corresponds to that
slope lies in the region between $B, T$ and, therefore, depends only
on the isomorphism class of $L/K$. So, our main theorem based on the previous lemmas can be described
below. Let $c_1, c_2$ be the constants obtained in the calculations in the previous section.

\begin{theorem}
Let $D(\mathbf{k},\alpha)$ be the number of eigenvalues of the $p^{- \sum
v_i}T_{p}$ operator acting on $\mathcal{S}^{D}_{\mathbf{k}}(U,R)$.
Let $ \alpha \leq c_1 n^{1/(d+1)} + c_2 \Rightarrow n \geq \lfloor
(\beta_{1} \alpha - \beta_{2})^{d+1} \rfloor = n(\alpha)$. Then, if
$\mathbf{k},\mathbf{k'} \geq n(\alpha)$, $\mathbf{k} \equiv \mathbf{k'} \mod p^{n(\alpha)}$ and
${\Gamma^{i}}$ are trivial, then $D(\mathbf{k'},\alpha) =
D(\mathbf{k'},\alpha)$. (The $\beta_i$ are constants which depend only on
$\alpha $ and $n$.)
\end{theorem}

\begin{proof}

We know from the previous lemma that if $\mathbf{k} \equiv \mathbf{k}' \mod p^{n-1}$, then
$\mathcal{S}^{D}_{\mathbf{k}}(\overline{U},R) / \mathcal{W}_{\mathbf{k}}^{D}(U,R)_{p}
\cong \mathcal{S}^{D}_{\mathbf{k}'}(\overline{U},R)/
\mathcal{W}_{\mathbf{k}'}^{D}(U,R)$. As $R$-modules we know that they are isomorphic to $\oplus_{i=1}^{r}
O/p^{a_{i}}O$. In the section on calculations, we have defined $B$
and $T$ such that $\alpha < c$ means that the number of eigenvalues
is the same (prop 3.2). Hence, we get that $D(\mathbf{k},\alpha) = D(\mathbf{k}',\alpha)$.

\end{proof}

$\mathbf{Note:}$ If the ${\Gamma^{i}}$ are not trivial, local constancy is
not possible. It is possible to get a result similar to the one we
get in the imaginary quadratic case i.e., an upper bound on
$D(\overrightarrow{k},\alpha)$ independent of the weight and depending only on the
slope. Yamagami has suggested that a technique of Buzzard in his paper titled "Eigenvarieties" can be used to overcome the obstructions.

\section{The Imaginary Quadratic case}

\subsection{Preliminaries}

Let $K$ be an imaginary quadratic field of class number $1$. Let $O$
be its ring of integers and $p$ be an odd rational prime in $K$,
which is inert in K. In Miyake's paper \cite{M}, one has a description of automorphic
forms for number fields $K$ in terms of $L_0^2$ decomposition and automorphic representations. We use Taylor's definitions \cite{T1}. Using the Eichler-Shimura isomorphism, one can translate this
description of automorphic forms to cohomology groups.

For any pair of non-negative integers $n_1, n_2$ we have a free
$(n_1+1)(n_2+1)$-dimensional $O$ module with an action of $GL_2(O)$
(or $M_2(O)$). It may be explicitly described as
$S_{n_{1}}(O^{2}) \otimes S_{n_{2}}(O^{2})$, where $S_{n}$ denotes
the $n$-th symmetric power , and where $\gamma \in GL_2(O)$ acts on
the first $O^{2}$ in the natural fashion and on the second via
$\overline{\gamma}$ (complex conjugation). We will denote this module $S_{n_{1}, n_{2}}$.
If $A$ is an $O$ module $S_{n_{1}, n_{2}}(A)$ will denote $S_{n_{1},
n_{2}}(O) \otimes (A)$. In particular, $S_{n_{1},
n_{2}}(\mathbb{C})$ is an irreducible finite dimensional
representation of the Lie group $SL_2(\mathbb{C})$.

We are interested in the cohomology of the congruence subgroups
$\Gamma < SL_2(O)$ with coefficients in $S_{n_{1}, n_{2}}(A)$. We let $\Gamma$ be a congruence subgroup of $SL_2(O)$ which is finitely generated. We refer the
reader to \cite{So} for the description.

We will consider $H^{1}(\Gamma, S_{n,n}(O))$ as our space of modular
forms over $K$. There is a notion of Hecke operators associated to
this space. Using Taylor's description they can be summarised as,

\begin{itemize}
  \item $[\Gamma_2 g \Gamma_1]: M^{\Gamma_1} \rightarrow
  M^{\Gamma_2}$ by $m \mapsto \sum (\gamma_i g)m$
  \item $[\Gamma_2 g \Gamma_1]: H^{1}(\Gamma_1, M) \rightarrow H^{1}(\Gamma_2,
  M)$ is induced by sending a $\Gamma_1$-cocycle $\phi$ to the
  $\Gamma_2$-cocycle $\delta \rightarrowtail \sum
(\gamma_{i}g)\phi((\gamma_{i}g)^{-1}\delta(\gamma_{j_{i}}g))$, where
$j_i$ is the unique index such that $\gamma_i^{-1} \delta
\gamma_{j_{i}} \in g \Gamma_i g^{-1}$
\end{itemize}

\subsection{Calculations}

Fix $n$ and let $g \geq n$.

Let $\Gamma = \Gamma_{1}(Np)$, and let $m$ be the minimal number of
generators of $\Gamma$. Let $S_{g,g} = S_{g}(O_{p}^{2})\otimes S_{g}(O_{p}^{2})$, which was described above in the previous section. We can think of $S_g$ in terms of homogeneous polynomials of degree
$g$. Let $O_p [x,y] = \oplus S_g (O_p^2)$. Let $J \subseteq O_p [x,y]$ denote the
homogeneous ideal $(p,x)$. Then for all $n > 0$, $J^{n}$ is also
homogeneous, so it can be written as $J^{n} = \oplus
M^{n}_{g}(O_p)$. The $M_{g}^{n}(O_p)$ are preserved under $\Gamma$.
Since
 if $\bp a &b\\c&d\ep \in \Gamma$ then $c \equiv 0 \mod p$.

We define $W_{g,g}^n = M_g^n \otimes S_{g}(O_{p}^2) + S_g (O_{p}^2)
\otimes M_g^n$

Let $I_{g,g}^n := S_g / M_g^n \otimes S_g / M_g^n
 \cong \oplus_{i=1}^{n^2} O/ p^{a_{i}}O$, where $1
\leq a_{i} \leq n$ and are arranged in decreasing order. Then, we have a short exact sequence,

$ 0 \rightarrow W_{g,g}^n \rightarrow S_{g,g} \rightarrow I_{g,g}^n
\rightarrow 0$.

This leads to a long exact sequence,

$H^{0}(\Gamma,I_{g,g}^n) \rightarrow H^{1}(\Gamma, W_{g,g}^n)
\rightarrow H^{1}(\Gamma, S_{g,g}) \rightarrow H^{1}(\Gamma,
I_{g,g}^n) \rightarrow H^{2}(\Gamma, W_{g,g}^n)$

Now, we take the maximal torsion-free quotient $TF$ and see that,

$0 \rightarrow H^{1}(\Gamma, W_{g,g}^n)^{TF} \rightarrow
H^{1}(\Gamma, S_{g,g})^{TF} \rightarrow H^{1}(\Gamma, I_{g,g}^n)^{*}
\rightarrow 0$, where $H^{1}(\Gamma, I_{g,g}^n)^{*}$ is a subquotient of $H^{1}(\Gamma, I_{g,g}^n)$

\smallskip

Consider, $L_{g} = H^{1}(\Gamma, S_{g,g})^{TF}$ and $K_{g}^n = H^{1}(\Gamma, W_{g,g}^n)^{TF}$.

Now, $L_{g}/K_{g}^n$ is a subquotient of $H^{1}(\Gamma, I_{g,g}^n)$,
which is itself a quotient of the group of $1$-cocycles from
$\Gamma$ to $I_{g,g}^n$. This group of cocycles is isomorphic to a
subgroup of $(I_{g,g}^{n})^m$, where $m$ is the minimal number of
generators of $\Gamma$ and $I_{g,g}^n \cong \oplus_{i=1}^{n^2} O/
p^{a_i}O$. Thus, $L_{g}/K_{g}^n \simeq \oplus_{i =1}^{s} (O/ p^{a_{i}^{g}}O)$,
where $1 \leq s \leq m n^2$

\smallskip
As in the case of Hilbert modular forms, we define $b_{i}^{g} = n - a_{i}^{g}$ and $B^{g}(j) = \sum_{i=1}^{j} b_{i}^{g}$.

\subsection{Lemmas}

\smallskip

Let $\xi = T_{p} = \Gamma \bp p &0\\0&1\ep \Gamma$. From \cite{T1} we know
$\Gamma \bp p &0\\0&1\ep \Gamma = \coprod \Gamma \alpha_{u}$, where
$\alpha_{u} = \bp p &u\\0&1\ep$ and $u$ runs over any set of
representatives for congruent class of $O$ mod $p$. If $\phi$ is a
$\Gamma$-cocycle, then it gets mapped to the cocycle $\delta \mapsto
\sum (\alpha_{u})\phi((\alpha_{u})^{-1}\delta \alpha_{v_u})$, where
$v_u$ is the unique $v$ such that $\alpha_u^{-1} \delta \alpha_{v_u} \in
\Gamma$.

\begin{lemma}
 $\xi(K_{g}^n) \subseteq p^{n}(L_{g})$.
\end{lemma}

\begin{proof}

Now, if $\kappa : \Gamma \rightarrow W_{g,g}^n $ is a $1$-cocycle,
it will be enough to show that $\xi(\kappa)$ is divisible by
$p^{g}$. By the definition of $\xi$, we need to show that $\bp p
&u\\0&1\ep M_{g}^n \subset p^{n}(S_{g})$. Let $f(x,y) \in M_g^n$. Then under the action of $\alpha = \bp p
&u\\0&1\ep$, $f(x,y) \mapsto f(ax + cy, bx + dy)$. Now, $f(x,y)$ is
a homogeneous polynomial of degree $g$, and $a,c \equiv 0 \mod p$ so
$f|\alpha = p^{g}f_1(x,y)$, where $f_1(x,y) \in S_{g,}$. Thus, we
get that $\xi(K_{g}^n) \subseteq p^{n}(L_{g})$.

\end{proof}

\smallskip

Let $\xi:L \rightarrow L$ be an $O_{p}$ linear endomorphism such
that $\xi(K) \subseteq p^{n}(L)$. Let $\overline{L}$ denote $L/K$
with its induced action of $\xi$, where $L/K \cong
\oplus_{i=1}^{r}O/p^{a_{i}}O$ where the $a_{i} \leq n$ and are
decreasing and define $b_i = n - a_{i}$ and $B(j) =
\sum_{i=1}^{j}b(i)$. Let $\sum_{i = 1}^t d_{i}X^{i}$ be the
characteristic polynomial of $\xi$ acting on L, and write $d_{s} =
0$ for $s>t$. From now on, we shall think of $L$ as an
$O_{p}[\xi]$ module, where $O_{p}[\xi]$ is thought of as an
indeterminate over $O_{p}$. So, $K$ is also a $O_{p}[\xi]$ submodule of
$L$.

Now, we state our key lemma.

\begin{lemma}
 $d_{i} \equiv 0 \mod p^{B(i)}$.
\end{lemma}

\begin{proof}
Choose an $O_{p}$ basis $(e_{i})$ of $L$ such that $(p^{a_{i}}
e_{i})$ is an $O_{p}$ basis for $M$. Let $(u_{i,j})$ be the matrix
for the action of $\xi$ on $L$ with respect to the basis we've
chosen. The assumption that $\xi(K) \subseteq p^{n}(L)$ means that
$(u_{i,j})p^{a_j}e_{j} = p^{n} w  $, where $w \in L$ implies that
$p^{b_{j}}$ divides $u_{i,j}$, since $b_j = n - a_j$.

By the definition of the characteristic polynomial $det (x -
u_{i,j}) = \sum d_s x^{t - s}$. Then, we can see that:

$d_s = (-1)^{s} \sum_{J \subseteq \{ 1,2,...,t\} \text{of size s}}
\sum_{\sigma \in Symm(J)} sgn(\sigma) \prod_{j \in J} u_{j, \sigma
(j)}$.

Since we know that $p^{b_{\sigma (j)}} $ divides $u_{j, \sigma (j)}$, we can rewrite $d_s$ as:

$d_s = (-1)^{s} \sum_{J \subseteq \{ 1,2,...,t\} \text {of size s}}
p^{\sum_{j \in J} b_j} \sum_{\sigma \in Symm(J)} sgn(\sigma)
\prod_{j \in J} u_{j, \sigma(j)}/p^{b_{\sigma (j)}}$.

Since, the $b_i$ are increasing and for all $J$ of size $s$,
$\sum_{j \in J} b_j \geq b_1 + b_2 + ... + b_s $. Hence, $d_s$ is divisible by $p^{B(s)}$ for all $s$ and the 
lemma is proved.

\end{proof}

\begin{lemma}

Let $P$ be a finite $p$ group and let $p^n P = 0$. If $Q$ is a
subquotient of $P$ then, for $0 \leq \mu \leq n - 1$,

$|p^{\mu}Q/p^{\mu + 1}Q| \leq |p^{\mu}P/p^{\mu + 1}P|$

\end{lemma}

Let $P = Z^1(\Gamma, I_{g}^n)=(I_{g}^n)^m$. The above lemma
says that if $L_g/K_g^n$ is a subquotient of $P$ (a finite $p$
group), then $B^g
> B$, where $B$ is the function defined for $P$,
and $B^g$ is defined for $L_g/K_g^n$. (We defined $B$ in the
previous section). This gives us the following result:

\begin{lemma}

 The Newton polygon of $\xi$ on $L_g$ has a uniform lower bound.
 This lower bound is a piecewise linear function,
which has slope $0$ for $0 \leq x \leq m$, slope $1$ for $m \leq x
\leq 4m$, and in general slope $r$ for $r^2 m \leq x \leq (r+1)^2
m $, $r<n$, and slope $n$ for $x \geq n^2 m$

\end{lemma}

Let $D(g,\alpha)$ be the number of eigenvalues of $T_p$ acting on
the space of automorphic forms of weight $g$, slope $\alpha$. We have the following result:

\begin{theorem}

$D(g,\alpha)$ has a uniform upper bound and is always less than
$[3m(\alpha+1)^2/2] m $.

\end{theorem}

\begin{proof}

To get an upper bound on $D(g,\alpha)$, we consider the Newton
polygon of the $T_p$ operator on the space of modular forms of
weight $g$. Since we know that the eigenvalues of slope $\alpha$
is given by the length of the projection on the $x$ axis of that
side, we find an upper bound for that projection.

 By the previous lemma, we know that the Newton polygon of
this operator is bounded below by the function $B$ which has slope
$r$ for $r^2 m \leq x \leq (r+1)^2 m $. We can bound $B$ from
below by the function $(2/3)(x/m)^{3/2} - x$ (see the previous
section for an explanation of how we get this lower bound). This
lower bound for the Newton polygon transforms into an upper bound
for the slope $\alpha$ subspaces, which proves the theorem.

\end{proof}

\section{Concluding Remarks}

We remark that using Jacquet-Langlands (see \cite{Hi2} for a
detailed exposition) we can translate the above theorems to results
on cuspidal Hilbert modular forms. Our $T_p$ operator has an extra
power of $p$ due to the determinant. In the imaginary quadratic
case, the torsion in the cohomology groups is an obstruction towards
getting local constancy of the slope $\alpha$ subspaces. The author
is in the process of making some computations in the classical case
to see the amount of torsion involved in the cohomology
groups, which will hopefully lead to some insight in the imaginary
quadratic case. For some examples of Hilbert Modular forms over real
quadratic fields see Dembele's papers as well as an ongoing project
to develop a database of such examples. Finally, similar results on local constancy for Hilbert modular forms
have been proved by Yamagami \cite{Y}, who kindly read through a draft and offered many useful suggestions.

Aftab Pande.

Dept of Mathematics, Brandeis University, Waltham, MA, 02454.

Email: aftab.pande@gmail.com

\end{document}